\documentclass[reqno,12pt]{amsart}
\usepackage{amsmath, amssymb, amsthm, amsfonts} 
\usepackage[english]{babel}
\usepackage{bbm}
\usepackage{graphicx}
\usepackage{url}
\usepackage{soul}
\usepackage{epstopdf}
\usepackage[ruled,vlined]{algorithm2e}
\usepackage[a4paper,bindingoffset=0.5cm,left=2cm,right=2cm,top=2.5cm,bottom=2cm,footskip=.8cm]{geometry}
\usepackage{rotating}
\usepackage{amsbsy,enumerate}
\usepackage{comment}
\usepackage{mathrsfs} 
\usepackage{enumitem}
\usepackage{cancel}
\usepackage{subfigure}\usepackage{booktabs}
\usepackage{siunitx}

\newcommand{\bs}{\boldsymbol}

\newcommand{\vb}{\vspace{3.2mm}}
\renewcommand{\hat}{\widehat}

\newcommand{\vertiii}[1]{{\left\vert\kern-0.25ex\left\vert\kern-0.25ex\left\vert #1 \right\vert\kern-0.25ex\right\vert\kern-0.25ex\right\vert}}
\allowdisplaybreaks
\usepackage{hyperref}

\allowdisplaybreaks

\usepackage[utf8]{inputenc}
\usepackage{type1cm}         
\usepackage{multicol}        
\usepackage[bottom]{footmisc}

\usepackage{newtxtext}       %
\usepackage{newtxmath}       

\usepackage{pgfplots}

\usepackage{tikz}
\usetikzlibrary{automata,positioning}
\usetikzlibrary{decorations.pathreplacing,decorations.markings,shapes.geometric}
\usetikzlibrary{shapes,arrows}
\usetikzlibrary{backgrounds,calc,positioning}

\usepackage{pgfplots}
\pgfplotsset{compat=1.9}

\usetikzlibrary{chains,shapes.multipart}
\usetikzlibrary{shapes,calc}
\usetikzlibrary{automata,positioning}

\usepackage{tikz}
\usetikzlibrary{automata,positioning}
\usetikzlibrary{decorations.pathreplacing,decorations.markings,shapes.geometric}
\usetikzlibrary{shapes,arrows}
\usetikzlibrary{backgrounds,calc,positioning}

\begin{document}

	\title[A maximum entropy-based estimator for dynamic random graphs]{Analysis of a maximum-entropy based estimator \\for dynamic random graph models}
\author[{D. Garlaschelli, \:M. Mandjes, \:F.P. Pijpers {\tiny \:and\:} J. Wang}]{Diego Garlaschelli, Michel Mandjes, Frank P. Pijpers {\tiny and} Jiesen Wang}

	\begin{abstract}
		We study dynamic random graphs in which the set of nodes is fixed, but edges evolve over time according to an underlying stochastic mechanism. Using a maximum-entropy approach, we define a probability distribution on graph trajectories that is consistent with observed constraints, capturing the inherent uncertainty in partially observed networks. We introduce a moment-based estimator for the parameters of this distribution and establish its statistical properties, such as consistency and asymptotic normality, with explicit formulas for the covariance structure. Numerical experiments demonstrate the estimator’s accuracy and robustness across various dynamic network scenarios. Our framework bridges probabilistic modeling and statistical inference in time-varying networks, providing practical tools for understanding and predicting complex edge dynamics.
  
\vb


\vb

\noindent
{\sc Affiliations.} 
DG is with IMT School for Advanced Studies, Piazza San Francesco 19, 55100 Lucca, Italy, and Lorentz Institute for Theoretical Physics, Leiden University, Niels Bohrweg 2, 2333 CA Leiden, The Netherlands.

\noindent MM is  with Mathematical Institute, Leiden University, Einsteinweg 55, 2333 CC Leiden, 
The Netherlands; he is also affiliated with (a)~Korteweg-de Vries Institute for Mathematics, University of Amsterdam, Amsterdam, The Netherlands, (b)~E{\sc urandom}, Eindhoven University of Technology, Eindhoven, The Netherlands, (c)~Amsterdam Business School, Faculty of Economics and Business, University of Amsterdam, Amsterdam, The Netherlands.  

\noindent FPP is with Korteweg-de Vries Institute for Mathematics, University of Amsterdam, Science Park 105-107, 1098 XG Amsterdam, The Netherlands, and Statistics Netherlands (CBS), Henri Faasdreef 312, 2492 JP The Hague, The Netherlands

\noindent JW is with Korteweg-de Vries Institute for Mathematics, University of Amsterdam, Science Park 105-107, 1098 XG Amsterdam, The Netherlands.

\medskip

\noindent Date: {\it \today}.

\vb

\noindent
{\sc Acknowledgments.} 
DG's research is supported by the European Union - NextGenerationEU - National Recovery and Resilience Plan (Piano Nazionale di Ripresa e Resilienza, PNRR), projects `SoBigData.it - Strengthening the Italian RI for Social Mining and Big Data Analytics' - Grant IR0000013 (n. 3264, 28/12/2021) and ``Reconstruction, Resilience and Recovery of Socio-Economic Networks'' RECON-NET EP\_FAIR\_005 - PE0000013 ``FAIR'' - PNRR M4C2 Investment 1.3.
MM's and JW's research has been funded supported by the European Union’s Horizon 2020 research and innovation programme under the Marie Sklodowska-Curie grant agreement no.\ 945045, and by the NWO Gravitation project {\tiny NETWORKS} under grant agreement no.\ 024.002.003. \includegraphics[height=1em]{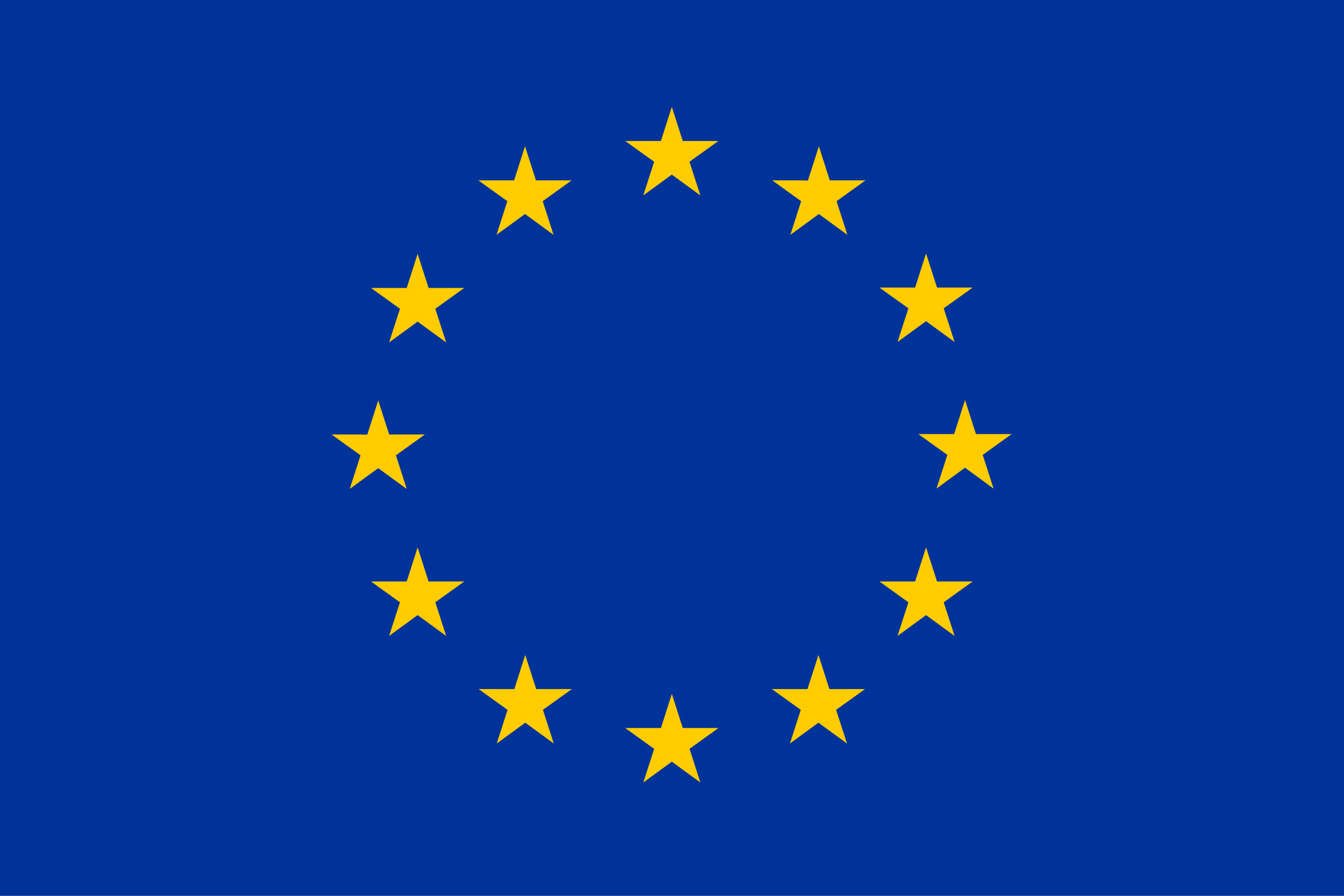}

\vb

\noindent
{\sc Email.} \url{diego.garlaschelli@imtlucca.it}, \url{m.r.h.mandjes@math.leidenuniv.nl}, \: 

\noindent
\url{f.p.pijpers@uva.nl},\:\url{j.wang2@uva.nl}

	\end{abstract}

 \maketitle

 \newpage

 \section{Introduction}
Where random graphs have traditionally been studied in a static context, there is an increasing realization that, in a wide range of application domains, the structure of the graph varies over time. In many real-world systems, both the set of nodes and the connections between them are inherently dynamic, driven by underlying stochastic mechanisms. Applications across social networks, communication systems, and biological networks highlight that interactions evolve on multiple time scales and may exhibit complex temporal dependencies. For example, the study of communication in highly dynamic networks, such as broadcasting and routing in delay-tolerant networks, requires models that capture temporal variability \cite{Vasilakos16}; similarly, opportunistic use of transportation networks exploits the underlying temporal patterns in mobility \cite{Pentland04}. In social and collaboration networks, dynamic interactions have been analyzed to reveal evolving community structures and coauthorship patterns \cite{Ghoshal09, Newman01, Newman04}. This shift in perspective has led to a surge of research on probabilistic models of time-varying graphs, which aim to faithfully represent the temporal nature of interactions while introducing new challenges in analysis, inference, and control \cite{Barabasi16, Durrett10, HolmeSaramaki12, PorterGleeson16}.

A central challenge in the study of dynamic random graphs is the characterization of their stochastic evolution. Given a time series of observations that summarize specific aspects of the graph, the goal is to infer the parameters governing the underlying probabilistic model. For example, one may observe, at discrete time points, the number of neighbors of each node, the presence or absence of particular edges, or other network summary statistics. Estimating the model parameters from such partially observed data is complicated by dependencies across nodes and across time, as well as by the fact that only aggregated information may be available.

The objective of this work is to develop an estimator for these parameters that not only is practically computable but also comes with rigorous statistical guarantees. A main objective of this paper is to demonstrate {\it consistency}:  the estimator that we propose converges to the true parameter as the number of observations, throughout this paper denoted by $T$, increases. Furthermore, we show that the estimator is {\it asymptotically normal}, meaning that for large $T$ the difference between the estimator and the true value behaves like a zero-mean normally distributed random variable, scaled by a factor of $1/\sqrt{T}$.

 \medskip

 This paper bridges two complementary streams of research in the study of dynamic networks: 
(a) the identification of model parameters using a {\it maximum-entropy} approach, which allows for the reconstruction of network structures from partial or aggregated observations in a principled, unbiased manner; and 
(b) the inference of these parameters using {\it statistical techniques}, including consistency analysis, asymptotic normality, and related probabilistic tools to quantify estimation uncertainty and construct confidence intervals. 

In the following, we provide a brief, non-exhaustive overview of the key ideas underlying both strands of research. The discussion is deliberately selective: we focus on concepts that are directly relevant to the methodology developed in this paper, while omitting related work that, though important, lies outside the immediate scope of our approach.

 \begin{itemize}
     \item[(a)]We begin by reviewing reconstruction methods for networks based on maximum entropy. The concept was first explored by Jaynes \cite{JAY}, who proposed inferring parameters using Shannon entropy by extending the statistical mechanics framework introduced by Gibbs. The central idea is to maximize the system's Shannon entropy over the unknown parameters, incorporating all available empirical information through the introduction of a set of constraints. Informally, this procedure `maximizes our ignorance about the system', in that it only uses the information that we are given through the observations. 
Over the past decades this maximum-entropy approach has developed into a powerful tool for assessing complex networks; for more background, with emphasis on economic and financial networks, we refer to the monograph \cite{GS}. 

\noindent The approach followed can be summarized as follows. Given a set of $N$ vertices, we define by $p(g):={\mathbb P}(G=g)$ the probability that the network is given by the graph $g\in\Omega$, where $\Omega$ is the set of all possible graphs on $N$ vertices.  Then the {\it Shannon entropy} is given by
\[-\sum_{g\in\Omega}p(g)\log p(g). \]
Suppose now that we are given information about the graph. Let ${\bs C}(g)$ a vector (of dimension $M\in{\mathbb N}$) of `properties' of the graph $g$, for instance the $N$-dimensional vector recording the number of neighbors of each of the vertices, and let $ {\bs c}$ its observed value. Then the maximum-entropy approaches solves, with the vector ${\bs\theta}\in{\mathbb R}^M$ representing the Lagrange multipliers pertaining to each of the properties, 
\[\max_{\bs\theta,\:{\bs p}(g)}\left(-\sum_{g\in\Omega}p(g)\log p(g) + \sum_{m=1}^M\theta_m \left(\sum_{g\in\Omega}p(g) \,C_m(g)- \bar c_m\right) \right).\]
It is directly verified that the optimizing probability distribution is given by
\[p^\star(g) ={\mathscr H}(g\,|\,{\bs\theta}):=\frac{1}{Z({\bs\theta)}}\exp(-H(g\,|\,{\bs\theta})), \]
with $H(g\,|\,{\bs\theta}):=\sum_m \theta_m   c_m$ denoting the Hamiltonian and $Z(\bs\theta)$ a normalizing constant. This distribution aligns with the functional form of exponential random graphs. 

\noindent  An example of a study in which this approach has been used, is \cite{CIM}. There the constraints are (weighted versions of) the per-vertex numbers of neighbors. The input of the procedure, i.e., the observed vector ${\bs c}$, is a snapshot of these quantities. The output is the probabilistic object $p^\star(g)$ that is such that sampling from it leads to network structures that reproduce the essentials of a network structure that `maximally align' with the observations. In the recent works \cite{GC, GC2} the maximum-entropy approach has been extended to the setting of dynamic graphs. 

\item[(b)]The second approach frames parameter identification as a statistical inference problem. As previously mentioned, the observations usually consist of a time series that captures specific characteristics of the dynamic graph, i.e., the detailed edge dynamics are not observed. For example, the data may represent periodic observations of the total number of edges in the graph. This makes our inference problem an {\it inverse problem}, where the goal is to infer the dynamics of the edges from observations of more aggregate processes.

\noindent
Inverse problems have been extensively studied across various classes of stochastic models. For example, in the context of branching processes, we refer to the textbook \cite{GUT}. Here, the challenge typically involves inferring characteristics of the offspring distribution from periodic observations of the population size. For inverse problems in queueing theory, we refer to the survey \cite{ASA}. A specific study in this domain is \cite{HAN}, which focuses on estimating the service-time distribution from workload observations in the M/G/1 queue. This is complicated by the fact that the workload process exhibits a complex dependence structure. Similarly, \cite{RAV} investigates a L\'evy-driven storage system, where the level of the system is observed at Poisson epochs, and the underlying L\'evy exponent is inferred. Recently, the topic of inference for dynamic graphs has garnered increasing interest. For example, the study \cite{MW} examines a dynamic version of the Erd\H{o}s-R\'enyi random graph, focusing on estimating the parameters related to the on- and off-time distributions of individual edges, based solely on subgraph counts.

\noindent A distinctive feature of these studies is the objective of rigorously establishing that the proposed estimator exhibits desirable statistical properties as the number of observations $T$ increases. One fundamental property is {\it consistency}, which ensures that the estimator converges to the true parameter vector as more data becomes available. Informally, this implies that with a sufficiently large amount of information, the parameters can be estimated with arbitrarily high accuracy. Beyond consistency, many estimators also admit {\it asymptotic normality}: in this case, the estimator can be expressed as the true parameter value plus a zero-mean normally distributed perturbation whose magnitude scales inversely with $\sqrt{T}$. Intuitively, this means that for large $T$, the estimator behaves like a random variable that is tightly concentrated around the true value, with fluctuations that decrease as more data is collected. Together, these properties provide both theoretical guarantees and practical guidance for inference in dynamic network models, giving confidence that the estimated parameters reliably reflect the underlying stochastic process.
 \end{itemize}

\medskip

This paper seeks to bridge the two complementary research areas. In line with branch (a), we consider dynamic graphs represented by a probability distribution obtained by maximizing the associated Shannon entropy, subject to constraints that encode the observed information about the network. This allows us to construct a principled probabilistic model of the evolving network, consistent with available data while remaining maximally non-committal about unobserved structure. In line with branch (b), our goal is to rigorously analyze the statistical properties of the resulting estimator, providing explicit performance guarantees. In particular, we establish asymptotic normality, showing that as the number of observations grows, the estimator converges to the true parameters and the estimation errors behave like a zero-mean Gaussian, appropriately scaled by the square root of the sample size $T$.

\medskip

This paper is organized as follows. In Section~\ref{S2}, we introduce our framework, detailing the dynamic random graph process, the observation scheme, and the construction of a probability distribution over the space of graph trajectories. Section~\ref{S3} then defines our estimator, which is essentially based on the method of moments, and rigorously establishes its statistical properties, including consistency and asymptotic normality; an explicit formula for the corresponding covariance matrix is also provided to facilitate practical implementation. Finally, in Section~\ref{S4}, we present numerical examples that illustrate the performance and robustness of the proposed estimator, highlighting its applicability to a variety of dynamic network scenarios.

 \section{Framework and preliminaries} \label{S2}
 In this section we subsequently introduce the notation that defines our dynamic graph process, discuss our observation process, define the probability distribution of the dynamic graph process, and present preliminaries on the edge dynamics. 
 \subsection{Dynamic random graph.} We consider a random graph, consisting of $N$ vertices. The graph is dynamic in the sense that the edges are evolving over time.
 {Time is discrete and indexed by the integer $t=1,\ldots,T$, where $T$ denotes the total number of time steps over which the graph dynamics is observed.}
 Let $G_t$ {denote the random graph at time $t$ and $g_t$ its particular realization. The `random graph trajectory' is given via
 \[{\mathscr G}_T \equiv \{G_1,\ldots, G_T\}\] and the particular realized `graph trajectory' is 
 \[{\mathscr g}_T \equiv \{g_1,\ldots, g_T\}.\]
 }
 We denote the {random} graph at time $t$ by
 \[G_t \equiv \{{A}_{ij}(t), \:i,j=1,\ldots, N\},\]
 {where $A_{ij}(t)$ is a random variable with $A_{ij}(t)=1$ denoting that the edge between the vertices $i$ and $j$ exists, and $A_{ij}(t)=0$ else, with $t=1,\ldots,T$. The realized value of  $A_{ij}(t)$ is denoted as  $a_{ij}(t)$.} We let $\Omega$ be the space of graphs with $N$ vertices, so that  ${{\mathscr g}_T}\in \Omega^T$. In the sequel we consider {\it undirected graphs}: in our setup the random graph's adjacency matrix is symmetric, i.e., ${A}_{ij}(t)={A}_{ji}(t)$ for all $i,j=1,\ldots,N$ and $t=1,\ldots,T;$ it is noted, however, that the asymmetric case can be dealt with analogously. All edges  are assumed to evolve independently. 

\subsection{Observation process.}
We throughout assume that we have access to the following information: for all vertices $i=1,\ldots,N$ we observe the {realized} sequences
\[\left(\sum_{j\not = i}a_{ij}(t)\right)_{t=1,\ldots,T},\quad \left(\sum_{j\not = i}a_{ij}(t)\,a_{ij}(t+1)\right)_{t=1,\ldots,T-1},\]
which are interpreted as the (a)~the number of vertices that were neighbor of vertex $i$ at any given point of time, and (b)~the number of vertices that were neighbor of vertex $i$ over two consecutive time steps (`one-lagged persisting edges'), respectively. 

 \subsection{Probability distribution.} \label{SS: distr} The next step is to define a probability distribution on the space of graph trajectories, in the following way. We deploy a dynamic random graph process of which the underlying dynamics are Markovian. For any $N$-dimensional vectors ${\bs\alpha}$ and $\bs\beta$, we first introduce
 the {\it Hamiltonian}: for any ${\mathscr g}_T\in\Omega^T$,
 \[H({\mathscr g}_T\,|\,\bs\alpha,\bs\beta):= \sum_{i=1}^n\alpha_i \bar k_i[{\mathscr g}_T]+ \sum_{i=1}^n\beta_i \bar h_i[{\mathscr g}_T],\]
 where
\begin{align*}
    \bar k_i[{\mathscr g}_T]&:=\frac{1}{T} \sum_{t=1}^T\sum_{j\not = i}a_{ij}(t),\\
       \bar h_i[{\mathscr g}_T]&:=\frac{1}{T-1}\sum_{t=1}^{T-1} \sum_{j\not = i}a_{ij}(t)\,a_{ij}(t+1).
\end{align*}
 Then we define the probability distribution 
 \[{\mathbb P}({\mathscr G}_T={\mathscr g}_T)={\mathscr H}({\mathscr g}_T \,|\,\bs\alpha,\bs\beta) := \frac{1}{Z(\bs\alpha,\bs\beta)}\exp\left(-H({\mathscr g}_T\,|\,\bs\alpha,\bs\beta)  \right),\]
where
 $Z(\bs\alpha,\bs\beta)$ is the normalizing constant, often referred to as the {\it partition function}; as pointed out in the introduction, this distribution maximizes the Shannon entropy conditional on it being in line with the observations. It directly follows that
 \[Z(\bs\alpha,\bs\beta)=\sum_{{\mathscr g}_T\in\Omega^T}\exp\left(-H({\mathscr g}_T\,|\,\bs\alpha,\bs\beta)  \right).\]
Observe that we parametrized the $N(N-1)/2$ edge processes by $2N$ parameters (namely, the entries of $\bs\alpha$ and $\bs\beta$).

 Due to the linear nature of $H({\mathscr g}_T\,|\,\bs\alpha,\bs\beta)$ (as a function of $\bs\alpha$ and $\bs\beta$, that is), it can be interpreted as a Laplace transform (up to a multiplicative constant). This observation directly entails that $\log Z(\bs\alpha,\bs\beta)$ is strictly convex (in $\bs\alpha$ and $\bs\beta$). Observe that we can rewrite
 \[{\mathbb P}({\mathscr G}_T={\mathscr g}_T)  = e^{-W{({\mathscr g}_T\,|\,\bs\alpha,\bs\beta)}
},\quad W{({\mathscr g}_T\,|\,\bs\alpha,\bs\beta)}:=\sum_{i=1}^n\alpha_i \bar k_i[{\mathscr g}_T]+ \sum_{i=1}^n\beta_i \bar h_i[{\mathscr g}_T]-\log Z(\bs\alpha,\bs\beta).\]
 Note that $W{({\mathscr g}_T\,|\,\bs\alpha,\bs\beta)}$ is, as a function of $\bs\alpha$ and $\bs\beta$, is strictly concave, being a difference between a linear function and a strictly convex function. It thus follows that $W{({\mathscr g}_T\,|\,\bs\alpha,\bs\beta)}$ has a unique maximizing $(\bs\alpha,\bs\beta)$, hence the $2N$ first order conditions
 \[\bar k_i[{\mathscr g}_T] ={\mathbb E}\,\bar k_i[{\mathscr G}_T]= \frac{\displaystyle {\partial}Z(\bs\alpha,\bs\beta)/{\partial \alpha_i}}{Z(\bs\alpha,\bs\beta)},\quad \bar h_i[{\mathscr g}_T]={\mathbb E}\,\bar h_i[{\mathscr G}_T] = \frac{\displaystyle {\partial Z(\bs\alpha,\bs\beta)}/{\partial \beta_i}}{Z(\bs\alpha,\bs\beta)}\]
 have a unique solution; in both equations the first equality is to be interpreted as the moment equation to be solved, whereas the second indicates the way to compute the mean under consideration in terms of the partition function ${Z(\bs\alpha,\bs\beta)}$. 

\subsection{Underlying Markovian dynamics.} The next goal is to describe the Markov chain corresponding to the per-edge dynamics. We assume these edge processes are all in stationarity. Specializing to the edge between the vertices $i$ and $j$, we use the notation
\[p_{ij}:={\mathbb E}{A}_{ij}(t) = {\mathbb P}({A}_{ij}(t)=1),\quad q_{ij}:={\mathbb E}[{A}_{ij}(t)\,{A}_{ij}(t-1)]={\mathbb P}({A}_{ij}(t)={A}_{ij}(t-1)=1). \]
From these $p_{ij}$ and $q_{ij}$ we can compute the per-edge transition probabilities. For instance, it is readily verified that
\[{\mathbb P}({A}_{ij}(t)=0\,|\,{A}_{ij}(t-1)=0) = \frac{1-2p_{ij}+q_{ij}}{1-p_{ij}},\]
where it has been used that ${\mathbb P}({A}_{ij}(t)=0,{A}_{ij}(t-1)=0)$ can be rewritten as
\[1-{\mathbb P}({A}_{ij}(t)=1)-{\mathbb P}({A}_{ij}(t-1)=1)+{\mathbb P}({A}_{ij}(t)={A}_{ij}(t-1)=1).\]
Along the same lines, the full transition matrix, on the state space $\{0,1\}$, can be computed: for the edge between the vertices $i$ and $j$ the transition matrix $P_{ij}$ is given by
\[P_{ij}={\mathsf P}(p_{ij},q_{ij}),\quad\mbox{with}\:\:\:{\mathsf P}(p,q):=\left(\begin{array}{cc}{\displaystyle \frac{1-2p+q}{1-p}}&
{\displaystyle \frac{p-q}{1-p}}\\ \vspace{-3mm}&\\1- {\displaystyle \frac{q}{p}}&{\displaystyle \frac{q}{p}}
\end{array}\right)=: \left(\begin{array}{cc}
{\mathfrak p}_{00}(p,q)&{\mathfrak p}_{01}(p,q)\\ \vspace{-3mm}&\\{\mathfrak p}_{10}(p,q)&{\mathfrak p}_{11}(p,q)
\end{array}\right)\]
As pointed out in \cite[Appendix 1]{GC} and \cite[Appendix B--E]{GC2}, the probabilities $p_{ij}$ and $q_{ij}$ can be expressed in terms of the parameter vectors $\bs\alpha$ and $\bs\beta$; via an appropriate reparametrization the model defined through our Hamiltonian can be rephrased in terms of a combination of non-interacting one-dimensional Ising models. To this end we define $x_i:=\exp(-\alpha_i/T)$, $y_i:=\exp(-\beta_i/T)$. Also, $x_{ij}:=x_ix_j$, $y_{ij}:=y_iy_j$; observe that $x_{ij}$ is a function of $\alpha_i$ and $\alpha_j$, and that $y_{ij}$ is a function of $\beta_i$ and $\beta_j$. In addition,
\[\lambda_{ij}^\pm\equiv \lambda_{ij}^\pm(x_{ij},y_{ij}):=e^{J_{ij}}\left(\cosh B_{ij}\pm \sqrt{\sinh^2B_{ij}+e^{-4J_{ij}}}\right)\]where\[J_{ij}:=\frac{\log y_{ij}}{4},\:\:B_{ij}:=\frac{\log(x_{ij}y_{ij})}{2}.\]
Having introduced this notation, we can present the closed-form expressions for $p_{ij}\equiv p_{ij}(\bs\alpha,\bs\beta)$ and $q_{ij}\equiv q_{ij}(\bs\alpha,\bs\beta)$ found in \cite{GC}. Defining the map 
\[\Phi(x,y,\ell_+,\ell_-):= \frac{xy-1}{2\sqrt{4x+(xy-1)^2}}\frac{\ell_+^T-\ell_-^T}{\ell_+^T+\ell_-^T}+\frac{1}{2},\]
it turns out that we have 
\[p_{ij}\equiv p_{ij}(\bs\alpha,\bs\beta)=\Phi(x_{ij},y_{ij},\lambda_{ij}^+,\lambda_{ij}^-).\]
Likewise, with
\[\Psi(\ell_+,\ell_-):= \frac{\ell_{-}^{T}\ell_+^{T-1}+\ell_{+}^{T}\ell_-^{T-1}}{\ell_+^T+\ell_-^T},\]
we have
\[q_{ij}\equiv q_{ij}({\bs\alpha},\bs\beta)= p_{ij}^2+p_{ij}(1-p_{ij})\,\Psi(\lambda_{ij}^+,\lambda_{ij}^-).\]
Importantly, also the normalizing constant (i.e., partition function) can be evaluated. From \cite[Appendix C]{GC2}, we have the following expression in terms of $\lambda_{ij}^\pm$: with an empty product being defined as $1$,
\begin{equation}\label{defZ}Z({\bs\alpha},{\bs\beta})= \prod_{i=1}^N\prod_{j=1}^{i-1} \big((\lambda_{ij}^+)^T+(\lambda_{ij}^-)^T).\end{equation}
Expressions for 
 \[ \frac{\displaystyle {\partial}\ln Z(\bs\alpha,\bs\beta)}{\partial \alpha_i}=\frac{\displaystyle {\partial}Z(\bs\alpha,\bs\beta)/{\partial \alpha_i}}{Z(\bs\alpha,\bs\beta)}={\mathbb E}\,\bar k_i[{\mathscr G}_T],\quad \frac{\displaystyle {\partial}\ln Z(\bs\alpha,\bs\beta)}{\partial \beta_i}= \frac{\displaystyle {\partial Z(\bs\alpha,\bs\beta)}/{\partial \beta_i}}{Z(\bs\alpha,\bs\beta)}={\mathbb E}\,\bar h_i[{\mathscr G}_T]\]
 can be found from the closed-form expressions for $\lambda_{ij}^\pm$ in a routine manner using standard symbolic software; see also \cite[Appendix D]{GC2}. Alternatively, a highly accurate numerical proxy can be produced in the evident manner:\begin{align*}
     \frac{{\partial}\ln Z(\bs\alpha,\bs\beta)}{\partial \alpha_i} &\approx \frac{1}{h}\big(\ln Z(\bs\alpha+{\bs e}_ih,\bs\beta)-\ln Z(\bs\alpha,\bs\beta)\big),\\\frac{{\partial}\ln Z(\bs\alpha,\bs\beta)}{\partial \beta_i} &\approx \frac{1}{h}\big(\ln Z(\bs\alpha,\bs\beta+{\bs e}_ih)-\ln Z(\bs\alpha,\bs\beta)\big).   
 \end{align*}
with ${\bs e}_i$ the $i$-th $N$-dimensional unit vector (i.e., its $i$-th entry being one, and all other entries zero), and $h$ small. 

\section{Estimator, and its properties}  \label{S3}

In this section we propose our estimator, and prove that it provides performance guarantees: in the regime where $T$ grows large, the estimator of the $2N$-dimensional paramater vector $(\bs\alpha,\bs\beta)$ is asymptotically normal. We in addition identify the corresponding covariance matrix. 

\subsection{Estimator.}
The estimator $(\hat {\bs\alpha}_T,\hat{\bs\beta}_T)$ is defined as the solution to the $2N$ moment equations, in which we equate ${\mathbb E}\,\bar k_i[{\mathscr G}_T]$ and ${\mathbb E}\,\bar h_i[{\mathscr G}_T]$ to their empirical counterparts. Concretely, $\hat {\bs\alpha}_T$ and $\hat{\bs\beta}_T$ solve the following $2N$ equations in equally many unknowns:  for $i=1,\ldots,N$,
\begin{equation} \label{eq:mom}
\bar k_i[{\mathscr G}_T] ={\mathbb E}\,\bar k_i[{\mathscr G}_T\,|\,\bs\alpha,\bs\beta],\quad \bar h_i[{\mathscr G}_T] = {\mathbb E}\,\bar h_i[{\mathscr G}_T\,|\,\bs\alpha,\bs\beta],    
\end{equation}
where
\[ {\mathbb E}\,\bar k_i[{\mathscr G}_T\,|\,\bs\alpha,\bs\beta]= \sum_{j\not =i} p_{ij}(\bs\alpha,\bs\beta)=:P_i(\bs\alpha,\bs\beta),\quad  {\mathbb E}\,\bar h_i[{\mathscr G}_T\,|\,\bs\alpha,\bs\beta]= \sum_{j\not =i} q_{ij}(\bs\alpha,\bs\beta)=:Q_i(\bs\alpha,\bs\beta).\]
By the considerations presented at the end of \S\ref{SS: distr}, which relied on the strict concavity of $W(\bs\alpha,\bs\beta)$, these $2N$ equations have a unique solution.

The idea is now to first prove that the vector $({\bs k}[{\mathscr G}_T\,|\,\bs\alpha,\bs\beta], {\bs h}[{\mathscr G}_T\,|\,\bs\alpha,\bs\beta])$ is asymptotically normal as $T\to\infty$ (after centering by the vector of means and dividing by $\sqrt{T}$, that is). This is essentially a consequence of the standard central limit theorem for multivariate stationary processes, with the main property to be verified being that all normalized relevant variances and covariances are finite. In the second step, this asymptotic normality is then inherited by the estimator $(\hat {\bs\alpha}_T,\hat{\bs\beta}_T)$, as can be established by applying the celebrated {\it delta method} \cite{V2000,VdV}.

We proceed by giving a few considerations regarding the solution of the moment equations for large $T$. Observe that implicitly $\lambda_{ij}^\pm$, $p_{ij}$ and $q_{ij}$ are functions of the $\alpha_i/T$ and $\beta_i/T$ (as $x_i$ and $y_i$ are). In addition, recall the known property that, for any $v,w\in{\mathbb R}$ such that $v\not=w$,
\[\lim_{T\to\infty} \frac{v^T+w^T}{(\max\{v,w\})^T}=1.\]
Hence, by the definition of $Z({\bs\alpha},{\bs\beta})$ as given in 
\eqref{defZ},
for some function $\xi$ that does {\it not} depend on $T$,
\[Z({\bs\alpha},{\bs\beta}) \approx \big(\xi({\bs \alpha}/T,{\bs \beta}/T)\big)^T,\]
where `$\approx$' means that the ratio of the left-hand side and right-hand side converges to 1. It thus follows that, in the regime that $T\to\infty$,
\begin{align*}
    {\mathbb E}\,\bar k_i[{\mathscr G}_T\,|\,\bs\alpha,\bs\beta]&=\frac{\displaystyle {\partial}\ln Z(\bs\alpha,\bs\beta)}{\partial \alpha_i} = T 
    \,\frac{\partial \xi(\alpha/T,\beta/T)}{\partial \alpha_i} \\&= T \cdot\frac{1}{T}\cdot 
    \left.\frac{\partial \xi(x,y)}{\partial x_i}\right|_{{\bs x}={\bs\alpha}/T,{\bs y}={\bs\beta}/T}=\left.\frac{\partial \xi(x,y)}{\partial x_i}\right|_{{\bs x}={\bs\alpha}/T,{\bs y}={\bs\beta}/T},
\end{align*}
and likewise 
\[ {\mathbb E}\,\bar h_i[{\mathscr G}_T\,|\,\bs\alpha,\bs\beta]=\frac{\displaystyle {\partial}\ln Z(\bs\alpha,\bs\beta)}{\partial \beta_i}=\left.\frac{\partial \xi(x,y)}{\partial y_i}\right|_{{\bs x}={\bs\alpha}/T,{\bs y}={\bs\beta}/T}.\]
Now note that the sample means $\bar k_i[{\mathscr G}_T]$ and $\bar h_i[{\mathscr G}_T]$, appearing in \eqref{eq:mom}, converge to their ergodic limits as $T\to\infty$. Conclude that the solution of the moment equations (i.e., $(\hat {\bs\alpha}_T,\hat{\bs\beta}_T)$) is essentially linear:  as $T\to\infty$,
\[\frac{\hat {\bs\alpha}_T}{T}\to \bar {\bs\alpha},\quad \frac{\hat {\bs\beta}_T}{T}\to \bar {\bs\beta}\]
for some $(\bar {\bs\alpha},\bar{\bs\beta})$. Also, observe that, as $T\to\infty$,
\[\Phi(x,y,\ell_-,\ell_+)\to  \frac{xy-1}{2\sqrt{4x+(xy-1)^2}}+\frac{1}{2},\quad \Psi(\ell_+,\ell_-)\to\frac{\ell_-}{\ell_+}.\]
In view of these facts, recalling that we will be working in the setting that $T$ is large, in the sequel we thus let $\lambda_{ij}^\pm$, $p_{ij}$ and $q_{ij}$ be defined as above, but with $x_i:=\exp(-\bar\alpha_i)$ and $\smash{y_i:=\exp(-\bar\beta_i),}$ in which setting $p_{ij}$ and $q_{ij}$ do not depend on $T$.

\subsection{Asymptotic normality of the moments.} 
As mentioned, the asymptotic normality of the $2N$-dimensional vector $({\bs k}[{\mathscr G}_T\,|\,\bs\alpha,\bs\beta], {\bs h}[{\mathscr G}_T\,|\,\bs\alpha,\bs\beta])$,  in the regime that $T$ grows large, is a consequence of the central limit theorem for stationary processes.
The only condition that needs to be checked is that the covariances pertaining to pairs of entries of  $({\bs k}[{\mathscr G}_T\,|\,\bs\alpha,\bs\beta], {\bs h}[{\mathscr G}_T\,|\,\bs\alpha,\bs\beta])$ decay to zero inversely proportionally to $T$; the corresponding proportionality constants are then the entries of the covariance matrix appearing in the asymptotic normality statement.

A useful (and easily verifiable) identity in this context shows that ${\mathbb C}{\rm ov}({A}_{ij}(s),{A}_{ij}(t) )$ decays geometrically: 
\[{\mathbb C}{\rm ov}({A}_{ij}(s),{A}_{ij}(t) )=\gamma(p_{ij},q_{ij},|t-s|),\]where\[\gamma(p,q,t):= \frac{\zeta(p,q)}{(\eta(p,q))^2}\big(1-\eta(p,q)\big)^{t}=p(1-p)\big(1-\eta(p,q)\big)^{t},\]with\begin{align*}\eta(p,q)&:={\mathfrak p}_{01}(p,q)+{\mathfrak p}_{10}(p,q)=\frac{p-q}{p(1-p)},\\\zeta(p,q)&:={\mathfrak p}_{01}(p,q)\,{\mathfrak p}_{10}(p,q)=\frac{(p-q)^2}{p(1-p)}.\end{align*}
It is observed that this relation holds irrespective of $s$ being smaller or larger than $t$, due to the reversibility of the Markov chain governed by the transition matrix ${\mathsf P}(p_{ij},q_{ij}).$  Along the same lines, it can be found that, for $t>s$,
\begin{equation}\label{extrans}{\mathbb P}({A}_{ij}(t)=1\,|{A}_{ij}(s)=1) = p_{ij}+(1-p_{ij})\big(1-\eta(p_{ij},q_{ij})\big)^{t-s}.\end{equation}


(a)~To deal with ${\mathbb V}{\rm ar}\,\bar k_i[{\mathscr G}_T]$, we note that, as $T\to\infty$,
\begin{align*}
    T\,{\mathbb V}{\rm ar}\,\bar k_i[{\mathscr G}_T] &= \frac{1}{T} \sum_{j\not=i} {\mathbb V}{\rm ar}\,\left(\sum_{t=1}^T {A}_{ij}(t)\right)\\
    &=\sum_{j\not=i} \left({\mathbb V}{\rm ar}\,{A}_{ij}(1)+\frac{2}{T}\sum_{t=2}^T(T-t+1)\,{\mathbb C}{\rm ov}({A}_{ij}(1),{A}_{ij}(t))\right)\\
    &\to\sum_{j\not=i} \left(\gamma(p_{ij},q_{ij},0) +  2\sum_{t=2}^\infty \gamma(p_{ij},q_{ij},t-1)\right)\\
    &=\sum_{j\not=i} p_{ij}(1-p_{ij})\left(1 + 2\frac{1-\eta(p_{ij},q_{ij})}{\eta(p_{ij},q_{ij})}\right)=:v_i(k),
\end{align*}
where in the first equality the independence of the edge processes has been used, and in the last step the geometric shape of $\gamma(p,q,t)$ (in $t$). 

\noindent 
(b)~As it turns out, ${\mathbb V}{\rm ar}\,\bar h_i[{\mathscr G}_T]$ can be handled in a similar manner. Indeed, 
\begin{align*}
    (T-1)\,{\mathbb V}{\rm ar}\,&\bar h_i[{\mathscr G}_T] =\frac{1}{T-1}\sum_{j\not=i}{\mathbb V}{\rm ar}\left(\sum_{t=1}^{T-1}{A}_{ij}(t)\,{A}_{ij}(t+1)\right)\\&=\sum_{j\not=i} \left({\mathbb V}{\rm ar}({A}_{ij}(1)\,{A}_{ij}(2))+\frac{2}{T-1}\sum_{t=2}^{T-1}(T-t)\,{\mathbb C}{\rm ov}({A}_{ij}(1)\,{A}_{ij}(2),{A}_{ij}(t)\,{A}_{ij}(t+1))\right)\\
    &\to \sum_{j\not=i} \left({\mathbb V}{\rm ar}({A}_{ij}(1)\,{A}_{ij}(2))+2\sum_{t=2}^{\infty}{\mathbb C}{\rm ov}({A}_{ij}(1)\,{A}_{ij}(2),{A}_{ij}(t)\,{A}_{ij}(t+1))\right)
\end{align*}
Then observe that ${\mathbb V}{\rm ar}({A}_{ij}(1)\,{A}_{ij}(2))=q_{ij}(1-q_{ij})$, and, for any $t=2,3,\ldots$,
\begin{align*}
   {\mathbb C}{\rm ov}({A}_{ij}(1)\,{A}_{ij}(2)&,{A}_{ij}(t)\,{A}_{ij}(t+1)) = {\mathbb E}[{A}_{ij}(1)\,{A}_{ij}(2)\,{A}_{ij}(t)\,{A}_{ij}(t+1)] -q_{ij}^2\\&= q_{ij}\,{\mathbb P}({A}_{ij}(t-2)=1\,|\,{A}_{ij}(0)=1)\,{\mathbb P}({A}_{ij}(1)=1\,|\,{A}_{ij}(0)=1) -q_{ij}^2,
   \end{align*}
   which, by applying \eqref{extrans}, can be rewritten as
\[q_{ij}\big(p_{ij}+(1-p_{ij})(1-\eta (p_{ij},q_{ij}))^{t-2}\big)
\frac{q_{ij}}{p_{ij}}-q_{ij}^2
=\frac{q_{ij}^2}{p_{ij}}(1-p_{ij})\big(1-\eta (p_{ij},q_{ij})\big)^{t-2}.\]
After elementary algebra and again using the geometric series, we find, as $T\to\infty$,
\[T\,{\mathbb V}{\rm ar}\,\bar h_i[{\mathscr G}_T]\to \sum_{j\not=i} \left(q_{ij}(1-q_{ij})+2\frac{(1-p_{ij})^2q_{ij}^2}{p_{ij}-q_{ij}}\right)=:v_i(h).\]
(c)~Now consider
\begin{align*}
    T\,{\mathbb C}{\rm ov}(\bar k_i[{\mathscr G}_T],\bar h_i[{\mathscr G}_T])&=\frac{1}{T-1}\sum_{j\not=i}\sum_{j'\not =i}{\mathbb C}{\rm ov}\left(\sum_{s=1}^T{A}_{ij}(s),\sum_{t=1}^{T-1}{A}_{ij'}(t)\,{A}_{ij'}(t+1)\right).
\end{align*}
Observe that in the double sum over $j=j'$ there is only a contribution due to terms with $j=j'$, due to the independence of the edge processes. As a consequence, the expression in the previous display equals
\[\frac{1}{T-1}\sum_{j\not=i}\sum_{s=1}^T\sum_{t=1}^{T-1}{\mathbb C}{\rm ov}\left({A}_{ij}(s),{A}_{ij}(t)\,{A}_{ij}(t+1)\right).\]
Now distinguish between on one hand the case that $s\in\{t,t+1\}$, and on the other hand the case that it is not. The contribution due to the former case can be verified to equal $2\,q_{ij}(1-p_{ij})$, whereas the contribution (for a given $j\not=i$) of the latter case is
\[2\,q_{ij}(1-p_{ij})\frac{1-\eta(p_{ij},q_{ij})}{\eta(p_{ij},q_{ij})}.\]
Upon combining the above, we thus end up with, as $T\to\infty$,
\[ T\,{\mathbb C}{\rm ov}(\bar k_i[{\mathscr G}_T],\bar h_i[{\mathscr G}_T])\to \sum_{j\not=i}2\,q_{ij}(1-p_{ij})\left(1+\frac{1-\eta(p_{ij},q_{ij})}{\eta(p_{ij},q_{ij})}\right)=:v_i(k,h).\]
The three remaining cases, in which $i\not=j$, are easier to deal with due to the independence of the edge processes. The underlying calculations essentially follow from the cases (a), (b), (c). 

\noindent
(d)~For $i\not=j$,  as $T\to\infty$,
\begin{align*}
    T\,{\mathbb C}{\rm ov}(\bar k_i[{\mathscr G}_T],\bar k_j[{\mathscr G}_T]) \to  p_{ij}(1-p_{ij})\left(1 + 2\frac{1-\eta(p_{ij},q_{ij})} {\eta(p_{ij},q_{ij})}\right) \\
    = p_{ij}(1-p_{ij})\left(\frac{2-\eta(p_{ij},q_{ij})}{\eta(p_{ij},q_{ij})}\right)=:c_{ij}^{(kk)}.
    \end{align*}
(e)~For $i\not=j$, as $T\to\infty$,
\begin{align*}
    T\,{\mathbb C}{\rm ov}(\bar h_i[{\mathscr G}_T],\bar h_j[{\mathscr G}_T]) \to q_{ij}(1-q_{ij})+2\frac{(1-p_{ij})^2q_{ij}^2}{p_{ij}-q_{ij}}=:c_{ij}^{(hh)}.
    \end{align*}
(f)~For $i\not=j$, as $T\to\infty$,
\begin{align*}
 T\,{\mathbb C}{\rm ov}(\bar k_i[{\mathscr G}_T],\bar h_j[{\mathscr G}_T]) \to
2\,q_{ij}(1-p_{ij})\left(1+\frac{1-\eta(p_{ij},q_{ij})}{\eta(p_{ij},q_{ij})}\right) \\
=2\,q_{ij}(1-p_{ij}) \left(\frac{1}{\eta(p_{ij},q_{ij})}\right)=:c_{ij}^{(kh)}
.
    \end{align*}
Define $C(k)$ as the $N\times N$ matrix with entries $c_{ij}^{(kk)}$, where it is noticed that its diagonal consists of zeroes, since the $A_{ii}=0$ (no self-loops). Define $C(h)$ and $C(k,h)$ analogously. Now define, in self-evident block-matrix notation, the covariance matrix (of dimension $2N\times2N$)
\[\Sigma:= \left(\begin{array}{cc}C(k)&C(k,h)\\C(k,h)&C(h)\end{array}\right)+ \left(\begin{array}{cc}{\rm diag}\{{\bs v}(k)\}&{\rm diag}\{{\bs v}(k,h)\}\\{\rm diag}\{{\bs v}(k,h)\}&{\rm diag}\{{\bs v}(h)\}\end{array}\right).\]
We have thus found the following distributional convergence of the moment equations our estimator is based on:
\begin{equation}\label{eqnorm}\sqrt{T}\left(\begin{array}{c}\bar {\bs k}[{\mathscr G}_T]-{\bs P}(\bs\alpha,\bs\beta)\\\bar {\bs h}[{\mathscr G}_T]-{\bs Q}(\bs\alpha,\bs\beta)\end{array}\right)\stackrel{\rm d}{\to} {\rm N}({\bs 0}, \Sigma),\end{equation}
where ${\bs 0}$ is the $2N$-dimensional all-zeroes vector, and ${\mathrm N}({\bs\mu},\Sigma)$ denotes a Normally distributed random vector with mean ${\bs\mu}$ and covariance matrix $\Sigma$. 

\subsection{Asymptotic normality of the estimator} The asymptotic normality of the estimator vector $(\hat {\bs\alpha}_T,\hat{\bs\beta}_T)$ follows by the asymptotic normality of $({\bs k}[{\mathscr G}_T\,|\,\bs\alpha,\bs\beta], {\bs h}[{\mathscr G}_T\,|\,\bs\alpha,\bs\beta])$ in combination with elementary Taylor expansions. The underlying principle is formally backed by the so-called {\em delta method}; see e.g.\ \cite[Ch.~XIX]{V2000} and \cite[\S 3.9]{VdV}. 

The starting point is provided by the moment equations
\begin{equation}\label{mom}\left(\begin{array}{c}{\bs P}(\hat{\bs\alpha}_T,\hat{\bs\beta}_T)\\{\bs Q}(\hat{\bs\alpha}_T,\hat{\bs\beta}_T)\end{array}\right)= \left(\begin{array}{c}\bar {\bs k}[{\mathscr G}_T]\\\bar {\bs h}[{\mathscr G}_T]\end{array}\right).\end{equation}
First consider the left-hand side of \eqref{mom}. By a first-order Taylor expansion at $(\bs\alpha,\bs\beta)$, neglecting higher order terms, we obtain
\[\left(\begin{array}{c}{\bs P}({\bs\alpha},{\bs\beta})\\{\bs Q}({\bs\alpha},{\bs\beta})\end{array}\right)+V(\bs\alpha,\bs\beta)\left(\begin{array}{c}\hat{\bs\alpha}_T-\bs\alpha\\\hat{\bs\beta}_T-\bs\beta\end{array}\right),\quad V\equiv V(\bs\alpha,\bs\beta):=\left(\begin{array}{cc}{\displaystyle \frac{\partial {\bs P}}{\partial\bs\alpha}}&{\displaystyle \frac{\partial {\bs P}}{\partial\bs\beta}}\\
\vspace{-3mm}&\\
{\displaystyle \frac{\partial {\bs Q}}{\partial\bs\alpha}}&{\displaystyle \frac{\partial {\bs Q}}{\partial\bs\beta}}\end{array}\right).\]
Now consider the right-hand side of \eqref{mom}. By \eqref{eqnorm}, we can rewrite it as, again neglecting higher order terms,
\[\left(\begin{array}{c}{\bs P}({\bs\alpha},{\bs\beta})\\{\bs Q}({\bs\alpha},{\bs\beta})\end{array}\right)+ \frac{1}{\sqrt{T}}\,{\rm N}({\bs 0},\Sigma).\]
Upon equating the two above expansions, we can identify the covariance matrix pertaining to the asymptotic normality of the estimator: we obtain the distributional identity, in the regime that $T\to\infty$,
\[\left(\begin{array}{c}\hat{\bs\alpha}_T\\\hat{\bs\beta}_T\end{array}\right) = \left(\begin{array}{c}{\bs\alpha}\\{\bs\beta}\end{array}\right)+\frac{1}{\sqrt{T}}\, V^{-1}\,{\rm N}({\bs 0},\Sigma).\]
We thus conclude that, as $T\to\infty$,
\[\sqrt{T}\left(\begin{array}{c}\hat{\bs\alpha}_T-\bs\alpha\\\hat{\bs\beta}_T-\bs\beta\end{array}\right)\stackrel{\rm d}{\to}{\rm N}\left({\bs 0}, V^{-1} \Sigma \,(V^{-1})^\top\right).\]
It is further noted that $({\bs P}(\alpha,\beta), {\bs Q}(\alpha,\beta))$ is to be interpreted as the derivatives of the log-Laplace transform of $(\bar{\bs k}[{\mathscr G}_T],\bar{\bs h}[{\mathscr G}_T])$ to all entries of $(\bs\alpha,\bs\beta)$. This means that $V({\bs\alpha},{\bs\beta})$ corresponds to the second derivative, or the Hessian matrix, which equals the covariance matrix $\Sigma$. As a consequence, 
\[V^{-1} \Sigma \,(V^{-1})^\top= \Sigma^{-1} \Sigma \,(\Sigma^{-1})^\top=\Sigma^{-1},\] so that,
as $T\to\infty$,
\[\sqrt{T}\left(\begin{array}{c}\hat{\bs\alpha}_T-\bs\alpha\\\hat{\bs\beta}_T-\bs\beta\end{array}\right)\stackrel{\rm d}{\to}{\rm N}\left({\bs 0}, \Sigma^{-1}\right).\]

\section{Numerical illustrations}
 \label{S4}

To illustrate the method of moments, we consider a parameter setting consisting of 50 pairs $(\alpha_i, \beta_i)$, for $i = 1,\ldots,50$.
The parameter values are generated as follows. We first draw them from a bivariate normal distribution. Specifically, for each realization, we consider a random vector
\[
(\alpha_i, \beta_i)^\top \sim \mathcal{N}(\mu, \Sigma),
\]
where the mean vector is given by
\[
{\boldsymbol \mu} = (0.5, 0.5),
\]
and the covariance matrix is
\[
\Sigma =
\begin{pmatrix}
0.1 & 0.05 \\
0.05 & 0.1
\end{pmatrix}.
\]

We then generate 50 independent realizations $\{(\alpha_i, \beta_i)\}_{i=1}^{50}$ from this distribution and treat them as fixed throughout the simulation study. We denote
\[
\bar{\bs \alpha} = (\alpha_1,\ldots,\alpha_{50})^\top, \qquad
\bar{\bs \beta} = (\beta_1,\ldots,\beta_{50})^\top.
\]

The non-zero off-diagonal entries of $\Sigma$ induce positive correlation between the two components.

Based on these fixed parameter values, we conduct numerical experiments with 100 independent replications.

\begin{figure}[ht]
\centering

\subfigure[Histogram of pooled standardized errors]{
    \includegraphics[width=0.45\textwidth]{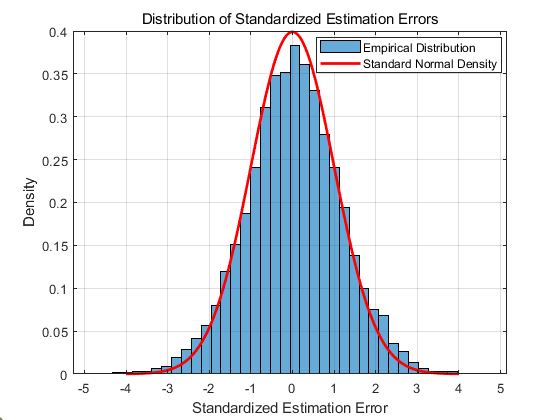}
}
\hfill
\subfigure[Q-Q plot against standard normal]{
    \includegraphics[width=0.45\textwidth]{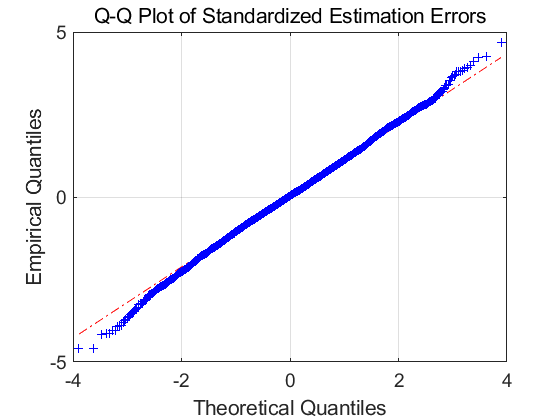}
}

\vspace{0.5cm}

\caption{
The histogram shows the empirical distribution of the pooled standardized estimation errors with the standard normal density superimposed.
The Q-Q plot compares the pooled standardized errors with the standard normal distribution.
}
\label{fig:diagnostics}

\end{figure}

In each replication, the estimator is constructed using observations over a time horizon of length $T = 10^5$. 
For each parameter component $i = 1,\ldots,50$ of $\smash{\bar{\bs\alpha}}$ and $\smash{\bar{\bs \beta}}$, and each replication $r = 1,\ldots,100$, we obtain estimates $\smash{\hat{\alpha}_i^{(r)}}$ and $\smash{\hat{\beta}_i^{(r)}}$. 
To allow meaningful comparisons across different parameter components, we standardize the estimation errors using the standard deviation:
\[
Z_i^{(r)} = \frac{\hat{\theta}_i^{(r)} - \theta_i}{\widehat{\sigma}_i},
\]
where $\theta_i$ denotes either $\bar {\bs\alpha}_i$ or $\bar {\bs\beta}_i$, and $\widehat{\sigma}_i$ is the sample standard deviation of the $100$ estimates for the $i$-th parameter. 
We then pool all standardized errors
\[
\left\{Z_i^{(r)} : i=1,\ldots,100, \ r=1,\ldots,100\right\}
\]
and analyze their empirical distribution. In particular, we present an empirical histogram against the standard normal distribution, as well as a Q-Q plot of the pooled standardized errors, to assess the validity of the asymptotic normal approximation. Figure \ref{fig:diagnostics} confirms the asymptotic normality. 

Furthermore, we summarize the overall estimation performance by reporting the mean and maximum root mean squared error (RMSE) across all 100 parameters; see Table \ref{tab:diagnostics}.

\begin{table}[ht]
\centering
\caption{Summary statistics of pooled standardized estimation errors and estimation accuracy.}
\label{tab:diagnostics}

\begin{tabular}{l S}
\toprule
Statistic & {Value} \\
\midrule
Mean of $Z$      & -0.001199 \\
Mean RMSE        & 0.006235 \\
Maximum RMSE     & 0.095221 \\
\bottomrule
\end{tabular}

\end{table}

 \end{document}